\documentclass[runningheads]{llncs}
\usepackage{graphicx}
\usepackage{amsmath,amssymb}
\usepackage{dsfont}
\usepackage{mathtools}
\usepackage{array}
\usepackage{multirow}
\usepackage{subfloat}
\usepackage{subcaption}
\usepackage{lscape}
\usepackage{booktabs}
\usepackage[hidelinks]{hyperref}
\usepackage[table]{xcolor}

\newcommand{\h}{\cellcolor{yellow!60!white}}
\newcolumntype{R}[1]{>{\raggedleft\let\newline\\\arraybackslash\hspace{0pt}}m{#1}}
\newcolumntype{L}[1]{>{\raggedright\let\newline\\\arraybackslash\hspace{0pt}}m{#1}}

\newcommand{\UBS}{\mathcal{U}}
\newcommand{\LBS}{\mathcal{L}}

\begin{document}
\title{Adapting Branching and Queuing for Multi-objective Branch and Bound}
%
%
\author{Julius Bauß \orcidID{0000-0003-2826-4369} \and
Michael Stiglmayr \orcidID{0000-0003-0926-1584}}
\authorrunning{J.~Bauß and M.~Stiglmayr}
%
\institute{Optimization Group, School of Mathematics and Natural Sciences, \\
University of Wuppertal, Wuppertal, Germany \\
\email{\{bauss,stiglmayr\}@uni-wuppertal.de}}
\maketitle              
\begin{abstract}
Branch and bound algorithms have to cope with several additional difficulties in the multi-objective case. Not only the bounding procedure is considerably weaker, but also the handling of upper and lower bound sets requires much more computational effort since both sets can be of exponential size. Thus, the order in which the subproblems are considered is of particular importance. Thereby, it is crucial not only to find efficient solutions as soon as possible but also to find a set of (efficient) solutions whose images are well distributed along the non-dominated frontier.
In this paper we evaluate the performance of multi-objective branch and bound algorithms depending on branching and queuing of subproblems. We use, e.g., the hypervolume indicator as a measure for the gap between lower and upper bound set to implement a multi-objective best-first strategy. We test our approaches on multi-objective knapsack and generalized assignment problems.

\end{abstract}

\section{Introduction}
In recent times multi-objective branch and bound algorithms have been investigated quite a lot for different types of optimization problems, e.g., integer, mixed-integer linear and mixed-integer nonlinear problems, see e.g.\ \cite{Belotti2013a,przybylski17multi}. Despite this increasing interest branch and bound algorithms suffer in the multi-objective setting from the weakness of bounding, the handling of large upper and lower bound sets and computational effort to compute them. Considering this, the algorithmic control of branch and bound is even more important. The creation and sequencing of subproblems can significantly increase the probability of finding efficient solutions in early stages of the algorithm. 

In this article we consider multi-objective integer optimization problems \eqref{eq:moilp} and their numerical solution using branch and bound algorithms. 
\begin{equation}\label{eq:moilp} \tag{MOILP}
   \begin{array}{rr@{\extracolsep{1ex}}c@{\extracolsep{1ex}}ll}
      \min \;\, & \multicolumn{3}{@{\extracolsep{0.75ex}}l}{\displaystyle \bigl( z_1(x),\ldots, z_p(x)\bigr)^\top}\\
      \mathrm{ s.t.}\;\, & \displaystyle A\,x &\leq&  b  \\
      &x &\geq& 0 \\
      &x &\in& \mathds{Z}^n 
   \end{array}
\end{equation}
where \(z(x) \coloneqq  (z_1(x),\ldots, z_p(x)) = C\cdot x\in\mathds{R}^p\) is the \(p\)-criterial linear objective function and $X\coloneqq \{x \in\mathds{Z}^n: A \leq b, x\geq 0 \}$ is the feasible set. Thereby, we rely on the Pareto concept of optimality which is based on the componentwise order. A feasible solution $x\in X$ is called \emph{efficient} if there is no other solution $\hat{x}\in X$ dominating it, i.e., $z(\hat{x})_i \leq z(x)_i$ for all \(i\in\{1,\ldots,p\}\) and at least one of these inequalities holds strict. 
The set of  efficient solutions is denoted by  $X_E$. By $Y_N= \{ z(x) \in Y\colon x\in X_E \}$ we denote the set of the non-dominated points in the objective space. As solution of a multi-objective optimization problem we consider the computation of a so-called \emph{minimal complete set}, i.e., the set of all non-dominated points and one (efficient) preimage for each of them. A self-contained introduction to multicriteria optimization is given in \cite{ehrgott05multicriteria}.


\section{Multi-objective Branch and Bound}\label{sec:bnb}
Branch and bound approaches subdivide problems recursively subdivided into subproblems. The recursion stops if a subproblem is discovered to be irrelevant for the computation of a minimal complete set. Subproblems are stored in a tree structure where the root node is associated with the original problem.
In each iteration one node is selected and its lower bound and upper bound are updated. The active node can be fathomed if the corresponding subproblem is either infeasible, solved or does not contain efficient solutions. Otherwise the corresponding problem is split into subproblems and the associated nodes are added as child nodes of the active node (branching). A comprehensive survey of multi-objective branch and bound algorithms is given in \cite{przybylski17multi}. 


Like in single-objective branch and bound algorithms lower bound sets are often determined by solving linear relaxations of the respective subproblem. In our branch and bound framework we rely on \emph{Benson's outer approximation algorithm} \cite{Benson1998an} to solve the resulting multi-objective linear problems. An alternative solution approach would be \emph{dichotomic search} \cite{Przybylski2010a}. In contrast to the single-objective case we obtain by the linear relaxation an lower bound which corresponds to its non-dominated set.

The so-called \emph{incumbent list} \(\UBS\) contains all integer feasible solutions obtained during the run of the algorithm which are not dominated by another feasible solution found so far. Candidates for the incumbent list are determined mainly by the extreme supported solutions of the lower bound sets which are checked for integer feasibility in each iteration. An integer feasible solution $\bar{x}\in X$ is then appended to the incumbent list, if there is no other solution \(x\in\UBS\) in it dominating \(\bar{x}\), i.\,e., $C(x)\leq C(\bar{x})$. 
If a new solution $\bar{x}$ is added to the incumbent list in turn all solutions in the incumbent list which are dominated by \(\bar{x}\) are removed. In this framework we start with an empty upper bound set. 




Having updated the upper and lower bound sets one can check if the respective subproblem can be fathomed by infeasibility, optimality or dominance. Thereby, infeasiblilty is only determined if the LP-relaxation is infeasible, as well. The subproblem is solved to optimality and thus fathomed, if lower bound set and upper bound set coincide. In the multi-objective case, however, this happens only if both upper and lower bound set consist of a single point. A node can be fathomed by dominance if no feasible solution in the respective subtree contributes to the minimal complete set, which is the case if all feasible solution in the subtree are (weakly) dominated by a solution in the incumbent list. A necessary condition is therefore, that all points $l \in \LBS $ in the lower bound set of the active node there is a point in the incumbent list $u\in \UBS $ such that $u \leqq l$.

\section{Sequencing of Subproblems}\label{sec:branch}
In the generic description of multi-objective branch and bound given in the previous section we omitted the sequencing of subproblems. The sequencing is based on the way the subproblems are generated (branching rule) and the ordering in which the generated subproblems are considered (node selection).

\subsection{Branching Rules}
Since we are considering purely integer optimization problems, a natural way to generate subproblems is to apply binary branching on one variable \(x_k\), i.e., create subproblems with feasible sets \(X'\coloneqq \{x \in\mathds{Z}^n: A \leq b, x_k<\hat{x}_k, x\geq 0 \}\) and \(X''\coloneqq \{x \in\mathds{Z}^n: A \leq b, x_k\geq \hat{x}_k, x\geq 0 \}\), where \(\hat{x}_k\) is determined based on the solution of the LP-relaxation. In the case of binary problems  (as in our numerical tests), this is equivalent to setting \(x_k=0\) and \(x_k=1\), respectively.


In the \emph{most often fractional rule} we count for every variable in how many extreme points of the lower bound set it attains fractional values. The variable which is most often fractional is chosen as branching variable.
In \cite{Belotti2013a} it is suggested to
sum up the distances to the next integer value in all extreme points of lower bound set for each variable. The variable with maximal total distance is used as branching variable in the \emph{how fractional rule}.

The branching rule \emph{sum of ratios} is proposed in \cite{Bazgan2009solving} for the multi-objective knapsack problem. Thereby, for every variable the ratios $c^k_i / w_i$ between objective function coefficients $c^k_i$ and the weight of the item $w_i$ 
are computed. The variable with the smallest sum of ratios is selected as branching variable. This branching scheme can be analogously applied on e.g.\ facility-location problems and generalized assignment problems, when opening costs and workload of a task are interpreted as the ``weight'' $w_i$. 
In \cite{jorge2010nouvelles} the same ratio vectors $(c^k_i / w_i)_{i=1,\ldots,p}$ are computed for all variables. The variable with the ratio vector which is least often dominated by other ratio vectors is selected as branching variable.

Note that the most of fractional and the how fractional rule are dynamic, i.e., depend on the current node, while the other two branching rules are static and depend only on properties of the problem instance.

\subsection{Node Selection}

Node selection strategies can be categorized in static and dynamic strategies. The most frequently applied static strategies are \textit{depth-first search} (LIFO) and \textit{breadth-first search} (FIFO). Both variants do not require additional computations and are easy to implement but do not adapt to the problem structure. 
Dynamic node selection strategies choose the active node depending on the gap between upper and lower bound. There are several approaches which mimic the \emph{best first} node selection strategy which is most frequently applied in single-objective branch and bound algorithms \cite{Morrison2016branch}. In all following methods the node with the largest gap is chosen. The methods mainly differ in the way the gap between the lower and upper bound set is measured.

The \emph{local hypervolume gap}, an approximated hypervolume measure, is suggested in \cite{Bauss2023augmenting} for node selection in bi-objective branch and bound, which can easily be extended to more criteria. In \cite{jesus2021design} the impact of the exact hypervolume gap on the performance of multi-objectve branch and bound is evaluated. However, the numerical results show that the evaluation of the hypervolume is computationally so demanding that its positive effects are compensated in terms of running time. To avoid this computational effort the hypervolume of of the largest search zone is approximated by the volume of the simplex spanned by the local upper bound \cite{Klamroth2015on} and the intersections of its dominance cone with the lower bound set.

Similarly, the largest \emph{hypervolume of a search-zone box} spanned by a local upper bound and the local ideal point of the lower bound set can be used to select the active node. The well known \emph{Hausdorff distance} between upper and lower bound set and the so-called
\emph{width of enclosure} \cite{eichfelder2021a} (i.e.\ the largest minimal distance in one objective function between upper and lower bound set) are alternative methods to measure the optimality gap and select an active node.

\section{Numerical Tests}\label{sec:num} 
We test and evaluate all combinations of node selection strategies and branching rules
presented in Section~\ref{sec:branch}, namely
the node selection strategies: depth-first (DF), breadth-first (BF), local hypervolume gap (HVG), volume of the largest search-zone box (HVB), Hausdorff distance of upper and lower bound set (HD), width of enclosure (WOE) and the branching rules most often fractional (MOF), how fractional (HF), sum of ratios (SR), dominance of ratios (DOM). 
Thereby, we consider two different problem classes: 3-objective knapsack problems (KP) with $n=30, 40, 50$ variables from \cite{Kirlik2014a} and randomly generated 3-objective generalized assignment problems (GAP) with $n=27,48,75$ variables, see \cite{Bauss2023gitinstances}.


The algorithms were implemented in Julia 1.9.0 and the linear relaxations (for the lower bound set) were solved with Bensolve 2.1 \cite{Loehne2017the}. The numerical test runs were executed on a single core of a 3.20 GHz Intel\textsuperscript{\textregistered}  Core\texttrademark\ i7-8700 CPU  with 32~GB RAM. The number of nodes and computation times are average values over 10 instances. 
%
%
In Table \ref{tab:knap} the numerical results for the 3-objective knapsack problem are shown. The numbers in brackets indicate how many (if not all) of the instances have been solved in the time limit of one hour. For all considered instance sizes the combination HVG-HF is the best choice w.r.t.\ the number of created branch and bound nodes. Due to the computation time of the local hypervolume gap, it does not lead to the best solution times. The combination DF-SR creates more nodes but is the best choice in terms of the total computation time for all considered instance sizes.
Table \ref{tab:gap} shows the numerical results of the corresponding generalized assignment problems. For the considered instances the combination HVG-MOF is the best choice regarding the number of created notes. Regarding the total computation time, the combination BF-MOF seems to be overall the best choice, although in the case $n=48$ the time is undercut by HVB-MOF. 
Note that all the favorable combinations for GAP use the most often fractional branching rule.

Overall, our numerical results indicate that the local hypervolume gap leads to a significant reduction in terms of created nodes and yields competitive running times.

\begin{table}[htbp!]
\scriptsize
\centering
\begin{tabular}{|L{1.8cm}|R{1.4cm}R{1.4cm}|R{1.4cm}R{1.4cm}|R{1.8cm}R{1.8cm}|}
\hline
\multirow{2}{*}{\begin{tabular}[l]{@{}l@{}}(KP)\end{tabular}} &
  \multicolumn{2}{c|}{$p=3,n=30$} &
  \multicolumn{2}{c|}{$p=3,n=40$} &
  \multicolumn{2}{c|}{$p=3,n=50$} \\ \cline{2-7} 
 &
  \multicolumn{1}{c}{nodes} &
  \multicolumn{1}{c|}{time(s)} &
  \multicolumn{1}{c}{nodes} &
  \multicolumn{1}{c|}{time (s)} &
  \multicolumn{1}{c}{nodes} &
  \multicolumn{1}{c|}{time (s)} \\ \hline
DF-MOF  & 23985.0  & 12.8004  & 138368.4  & 111.1310  & 391170.2 & 438.8650
 \\
DF-HF   & 25386.8  & 13.5571  & 145034.2  & 115.6625  & 390635.6 & 436.7868
 \\
DF-SR   & 13251.6  & \h 7.5527& 51868.0   &\h 41.9066 & 116711.2 &\h 135.2845
 \\
DF-DOM  & 13562.6  & 7.7311   & 55519.6   & 44.8569   & 123384.0 & 144.2920
 \\
BF-MOF  & 22689.6  & 17.7095  & 149388.0  & 196.9730  & 432790.0 & 757.8097
 \\
BF-HF   & 21375.0  & 16.2624  & 134515.2  & 179.8361  & 407798.4 & 739.3253
 \\
BF-SR   & 341674.8 & 200.2637 & 1269261.8 & 1256.6632 & 1433460.3 (5) & 3187.8977 (5)
 \\
BF-DOM  & 329000.8 & 192.8351 & 1234688.4 & 1222.4330 & 1525962.0 (4) & 3393.6130 (4)
\\
HVG-MOF & 10355.4  & 9.4816   & 49898.4   & 88.4784   & 112975.4 & 334.5662
 \\
HVG-HF  & \h 9886.0& 8.7751   &\h 49432.8 & 81.3003   &\h 109233.2 & 307.8950
 \\
HVG-SR  & 27234.4  & 28.2833  & 101070.4  & 215.4569  & 206470.0 & 913.4311
 \\
HVG-DOM & 24432.4  & 25.3734  & 93457.2   & 199.2275  & 196919.6 & 873.1798
 \\
HVB-MOF & 14355.6  & 11.1059  & 94993.8   & 146.8714  & 233845.0 & 483.6893
 \\
HVB-HF  & 14241.4  & 10.4059  & 91549.2   & 130.7214  & 230877.6 & 462.5850
 \\
HVB-SR  & 36433.6  & 32.3518  & 151624.6  & 229.3973  & 325188.8 (9) & 930.8374 (9)
 \\
HVB-DOM & 35701.8  & 31.7019  & 148932.2  & 225.3204  & 348987.4 (9) &998.9598 (9)
 \\
HD-MOF  & 12462.2  & 9.6071   & 79112.4   & 126.2463  & 221310.9 (9)&821.2436 (9)
 \\
HD-HF   & 11919.6  & 9.0139   & 78892.8   & 127.9482  & 241806.8 & 826.3442
 \\
HD-SR   & 38704.4  & 38.9071  & 191524.8  & 571.2338  & 352888.6 (7) &2144.3936 (7)
 \\
HD-DOM  & 33764.6  & 33.9414  & 162471.4  & 484.5804  & 335842.2 (7) &2040.8080 (7)
 \\
WOE-MOF & 15256.2  & 10.5954  & 82633.0   & 99.8712   & 238207.2 & 423.3763
 \\
WOE-HF  & 15068.8  & 10.2355  & 83082.8   & 101.3736  & 239090.6 & 434.7994
 \\
WOE-SR  & 37537.6  & 25.3568  & 264445.2  & 337.7509  & 721534.2 & 1628.2839
 \\
WOE-DOM & 28092.8  & 18.9768  & 184564.6  & 235.7279  & 563025.8 (9) &1270.5786 (9)
 \\ \hline
\end{tabular}\vspace*{4pt}
\caption{Tri-objective knapsack problems \cite{Kirlik2014a}.} \label{tab:knap}
\end{table}

\begin{table}[htbp!]
\scriptsize
\centering
\begin{tabular}{|L{1.8cm}|R{1.4cm}R{1.4cm}|R{1.4cm}R{1.4cm}|R{1.8cm}R{1.8cm}|}
\hline
\multirow{2}{*}{\begin{tabular}[l]{@{}l@{}}(GAP)\end{tabular}} &
  \multicolumn{2}{c|}{$p=3,n=27$} &
  \multicolumn{2}{c|}{$p=3,n=48$} &
  \multicolumn{2}{c|}{$p=3,n=75$} \\ \cline{2-7} 
 &
  \multicolumn{1}{c}{nodes} &
  \multicolumn{1}{c|}{time(s)} &
  \multicolumn{1}{c}{nodes} &
  \multicolumn{1}{c|}{time (s)} &
  \multicolumn{1}{c}{nodes} &
  \multicolumn{1}{c|}{time (s)} \\ \hline
DF-MOF        & 2145.2 & 0.8989   & 25214.0 & 22.7623  & 150039.4 & 243.4742 \\
DF-HF         & 2312.8 & 0.9342   & 25624.6 & 23.0392  & 148997.8 & 237.6341 \\
DF-SR         & 2728.8 & 1.0449   & 31214.0 & 25.7561  & 187934.2 & 306.8806 \\
DF-DOM        & 2814.2 & 1.0776   & 28019.0 & 23.1198  & 150022.8 & 227.3948 \\
BF-MOF        & 1794.2 &\h 0.7826   & 17324.0 & 18.0416  & 83273.4  &\h 162.8612 \\
BF-HF         & 1926.0 & 0.8348   & 17614.2 & 18.6211  & 85289.2  & 174.1150 \\
BF-SR         & 2472.4 & 1.0170   & 23418.4 & 22.9168  & 113958.0 & 230.9721 \\
BF-DOM        & 2617.4 & 1.0767   & 20730.6 & 20.2865  & 97610.8  & 188.5939 \\
HVG-MOF       &\h 1778.4 & 0.8929   &\h 17030.6 & 22.7282  &\h 80720.2  & 276.1176 \\
HVG-HF        & 1910.8 & 0.9563   & 17330.4 & 23.0759  & 82971.2  & 283.6493 \\
HVG-SR        & 2430.2 & 1.1571   & 23147.0 & 28.0784  & 111541.8 & 405.7946 \\
HVG-DOM       & 2585.2 & 1.2309   & 20442.2 & 24.7974  & 95689.8  & 320.4952 \\
HVB-MOF       & 1802.0 & 0.8076   & 17285.4 &\h 17.9048  & 83615.4  & 166.1450 \\
HVB-HF        & 1927.8 & 0.8504   & 17608.8 & 18.1974  & 85378.8  & 170.1354 \\
HVB-SR        & 2450.8 & 1.0291   & 23538.2 & 22.3729  & 113930.0 & 228.2173 \\
HVB-DOM       & 2609.6 & 1.0958   & 20879.0 & 19.8453  & 98906.0  & 188.6762 \\
HD-MOF        & 2014.2 & 0.8762   & 23257.2 & 22.6969  & 132736.8 & 241.9842 \\
HD-HF         & 2185.6 & 0.9402   & 23755.2 & 23.0600  & 129838.2 & 229.5523 \\
HD-SR         & 2583.0 & 1.0513   & 30504.8 & 27.5358  & 172742.2 & 306.2228 \\
HD-DOM        & 2716.8 & 1.1058   & 26249.4 & 23.6946  & 132324.4 & 223.6301 \\
WOE-MOF       & 1891.6 & 0.9198   & 18414.0 & 23.1183  & 91496.4  & 280.4083 \\
WOE-HF        & 2013.8 & 0.9715   & 18764.4 & 23.3688  & 91927.6  & 283.7661 \\
WOE-SR        & 2502.0 & 1.1550   & 24475.2 & 27.8896  & 123246.4 & 400.8106 \\
WOE-DOM       & 2627.6 & 1.2130   & 21986.2 & 25.0534  & 106589.6 & 314.8106
 \\ \hline
\end{tabular}\vspace*{4pt}
\caption{Tri-objective generalized assignment problems \cite{Bauss2023gitinstances}.}
\label{tab:gap}
\end{table}

\vspace{-4em}
\bibliographystyle{splncs04}
\bibliography{mybibliography}

\end{document}